%% file: main.tex
\documentclass[letterpaper, 10pt, conference, english]{ieeeconf}   % Comment this line out
\IEEEoverridecommandlockouts                              % This command is only
\usepackage[top=50mm, bottom=50mm, left=50mm, right=50mm]{geometry}
\usepackage{lineno}
\usepackage{amssymb}
\usepackage{amsmath}
\usepackage{epsfig}
\usepackage{graphicx}
\usepackage{graphics}
\usepackage{float}
\usepackage{subfigure}
\usepackage{multirow}
\usepackage{color}
\usepackage{lineno}
\usepackage{fullpage}
\usepackage[normalem]{ulem} 
\usepackage{makeidx}
\usepackage{xspace}
\usepackage{wrapfig}
\makeindex

\newtheorem{theorem}{Theorem}

\definecolor{darkred}{rgb}{1, 0.1, 0.3}
\definecolor{darkblue}{rgb}{0.1, 0.1, 1}
\definecolor{darkgreen}{rgb}{0,0.6,0.5}

\newcommand {\mm}[1] {\ifmmode{#1}\else{\mbox{\(#1\)}}\fi}

\title{Data-Driven Distributionally Robust Electric Vehicle Balancing for Mobility-on-Demand Systems under Demand and Supply Uncertainties}
 
\author{
{Sihong He} \and {Lynn Pepin} \and {Guang Wang} \and {Desheng Zhang} \and {Fei Miao}
\thanks{This work has been published in International Conference on Intelligent Robots and Systems (IROS 2020). Sihong~He, Lynn~Pepin, and Fei~Miao are with the Department of Computer Science and Engineering, University of Connecticut, Storrs Mansfield, CT, USA 06268.  Email: \{sihong.he, lynn.pepin, fei.miao\}@uconn.edu. This work is also partially supported by NSF SAS-1849238 and CPS-1932223. Guang Wang and Desheng Zhang are with the Department of Computer Science, Rutgers University, Piscataway, NJ, USA 08901. Email: \{gw255, desheng.zhang\}@cs.rutgers.edu.}
}

\begin{document}
\maketitle
%\linenumbers
% \setcounter{page}{0}

\thispagestyle{plain}
\pagestyle{plain}

\begin{abstract}
As electric vehicle (EV) technologies become mature, EV has been rapidly adopted in modern transportation systems, and is expected to provide future autonomous mobility-on-demand (AMoD) service with economic and societal benefits. However, EVs require frequent recharges due to their limited and unpredictable cruising ranges, and they have to be managed efficiently given the dynamic charging process. It is urgent and challenging to investigate a computationally efficient algorithm that provide EV AMoD system performance guarantees under model uncertainties, instead of using heuristic demand or charging models. To accomplish this goal, this work designs a data-driven distributionally robust optimization approach for vehicle supply-demand ratio and charging station utilization balancing, while minimizing the worst-case expected cost considering both passenger mobility demand uncertainties and EV supply uncertainties. We then derive an equivalent computationally tractable form for solving the distributionally robust problem in a computationally efficient way under ellipsoid uncertainty sets constructed from data. Based on E-taxi system data of Shenzhen city, we show that the average total balancing cost is reduced by 14.49\%, the  average  unfairness  of  supply-demand ratio and utilization is reduced by 15.78\% and 34.51\% respectively with the distributionally robust vehicle balancing method, compared with solutions which do not consider model uncertainties. %This is about a 60-million-mile or a 8-million-dollar cost reduction annually in NYC.
\end{abstract}

\input{intro}

\section{Problem Formulation}
\label{sec:prob_form}
In this section, we formulate a distributionally robust optimization problem to balance EVs across a city with minimum total idle distance and balanced charging station utilization. Both passenger mobility demand and EV supply uncertainties are considered. The region every empty EV will go is updated in a receding horizon control process. At each time step, the EV status is updated to the dispatch center first, then the dispatch center calculates a vehicle balancing decision by solving the proposed distributionally robust optimization problem, and sent solutions to EVs. The goal is to dispatch vacant EVs to different regions to pick up current and predicted passengers if the EVs have enough energy, or to charging stations if the EVs are short of energy, while minimize the cost of dispatching for the following $\tau$ time steps. Local dispatchers that match individual EV with one or several passengers (for carpool) is out the scope of this work.%, but for this work we only focus on the centralized system-level vehicle balancing algorithm design. 

\subsection{EV States and Corresponding Actions}
We assume there are three possible states for one EV: vacant, occupied, and low-battery. Vacant means there are no passengers in this EV, and it has enough energy to finish the next trip. The controller dispatches vacant EVs according to current and predicted passengers demand. When a vacant EV picks up one or more passengers, it turns to occupied, and the controller has no actions for it until it becomes vacant again. An occupied EV will be finishing current order in a time period and will become a vacant EV once it drops off its passengers. One occupied EV can only become a vacant EV when it finishes the current order. When a vacant EV can not finish the next trip with the remaining battery, this EV becomes a low-battery EV and will go to regions assigned by the controller where it can find a charging station. Before a low-battery EV gets fully charged, it stays in the low-battery status until it leaves the charging station and becomes vacant. A low-battery can only transfer to a vacant EV or stay in current state. 

\subsection{Problem Description}
We assume that one day is divided into $K$ time intervals, and we use $k= 1,2,...,K$ to denote time index. We assume the entire city is divided into $N$ regions and we use $n$ to denote region index, where $n = 1,2,...,N$. At time $k$, the system-level controller makes vacant and low-battery EVs to go to other regions or stay in the same region for picking up passengers or charging, respectively. After one EV arrives at its dispatched region, a local-level controller assign the EV to pick up passengers or to charge according to the EV's battery status.

During time $k$, there are $r_i^k$ predicted total amount of passengers demand and $c_i^k$ predicted total number of EVs finish charging (new supply of EVs) in region $i$, where $i = 1,2,...,N, k = 1,2,...,K$. Let demand vector $r^k = [r_1^k, r_2^k,...,r_N^k]^T$ and supply vector $c^k = [c_1^k, c_2^k,...,c_N^k]^T\in \mathbb{R}^N$ be random vectors instead of deterministic vectors. And assuming they are independent. To model the spatial and temporal relations of deman(supply) during every $\tau$ consecutive time interval, we define concatenation of demand as $r = [r^1, r^2, ... , r^{\tau}]$,
and concatenation of supply as $c = [c^1, c^2, ... , c^{\tau}]$.
We use $F^*_r$ and $F^*_c$ to denote the unknown true probability distributions of $r,c \in \mathbb{R}^{N\tau}$ respectively, i.e. $r \sim F^*_r$ and $c \sim F^*_c$. 

We use non-negative matrices $X^k$ and $Y^k$ as the decision matrices at time $k$ where $X^k, Y^k \in \mathbb{R^{N\times N}_{+}}$ and $x_{ij}^k(y_{ij}^k)$ is the total amount of vacant(low-battery) EVs will be dispatched from region $i$ to region $j$ at the beginning of time $k$. Minimizing the expected allocating cost given true probability distributions of demand vector and supply vector is defined as the following stochastic programming problem:
\begin{align}
	\begin{split}
		\underset{X^{1:\tau}, Y^{1:\tau}}{\text{min.}}\ &\mathbb{E}_{r\sim F^*_r,c\sim F^*_c}\left[J(X^{1:\tau}, Y^{1:\tau},r,c) \right] \\
		\text{s.t.}\quad  & X^{1:\tau}, Y^{1:\tau}\in \mathcal{D},
	\end{split}
	\label{stochastic}
\end{align}
where $J(X^{1:\tau}, Y^{1:\tau},r,c)$ is a cost function of allocating EVs according to decisions $X^{1:\tau} = \{X^1, X^2,...,X^{\tau}\}$ and $Y^{1:\tau} = \{Y^1, Y^2,...,Y^{\tau}\}$ under demand $r$, supply $c$ and convex constraints domain  of decision variables $\mathcal{D}$.

However, in real world scenarios, we usually have limited knowledge about the true probability distributions of $r$ and $c$. Though we have historical or streaming data, we can only estimate a set of probability distributions, such that $F^*_r \in \mathcal{F}_c$, $F^*_c \in \mathcal{F}_c$ considering the randomness of the parameters and prediction errors~\cite{Tian2016Real, wang2019shared, EVAMoD_tcns19}, instead of knowing the exactly form of $F^*_r$ and $F^*_c$. Meanwhile, we notice that problem \eqref{stochastic} is computationally expensive to solve. So in this work, we consider to minimize the worst case expected cost function which is a minmax form of problem \eqref{stochastic} as below: \eqref{stochastic} as below:
\begin{align}
	\begin{split}
		\underset{X^{1:\tau},Y^{1:\tau}}{\text{min.}}\ \underset{F_r\in \mathcal{F}_r,F_c\in \mathcal{F}_c}{\text{max.}}\quad &\mathbb{E} \left[J(X^{1:\tau}, Y^{1:\tau},r,c) \right]\\
		%=\sum_{k=1}^{\tau} (J_D(X^k)+\beta J_E(X^{[1,\tau]},r^k))\right]\\
\text{s.t.}\quad  &X^{1:\tau}, Y^{1:\tau}\in \mathcal{D}.
	\end{split}
	\label{minmax}
\end{align}
Problem \eqref{minmax} is a form of distributionally robust optimization problem~\cite{Ye_dro} and assume $F^*_r \in \mathcal{F}_c$, $F^*_c \in \mathcal{F}_c$. In the following sections, we will define the complete forms of object function and constraints as well as the probability uncertainty sets $\mathcal{F}_r, \mathcal{F}_c$.

\subsection{Cost of EV Balancing $J_D$} 
We aim to balance the EV supply according to passenger mobility demand, by sending vacant state EVs to serve passengers according to $X^k$ and low-battery EVs to charging stations according to $Y^k$. Given a specific region partition method, let $W \in \mathbb{R}^{N \times N}$ be the cost matrix where $w_{ij}$ is the cost sending a vacant EV from region $i$ to region $j$. The cost can be metrics such as the approximated distance, the minimal routing distance or travel time between two regions. Here, we use approximated distance from two regions as the dispatching cost. Let $W^* \in \mathbb{R}^{N \times N} $ be the cost matrix of sending one low-battery EV which is partially the same as $W$. When there is at least one charging station in region $j$, $w^*_{ij} = w_{ij}$ for all $i$. If there are no charging stations in region $j$, $w^*_{ij} = \infty$ for all $i$ since low-battery EV should not go to regions without charging stations. 

Then the total re-balancing cost function $J_D$ for $\tau$ intervals is
\begin{align}
	J_D (X^{1:\tau}, Y^{1:\tau}) = \sum\limits_{k=1}^{\tau} \sum\limits_{i=1}^{N} \sum\limits_{j=1}^{N} (x^k_{ij} w_{ij} + \beta y^k_{ij} w^*_{ij}),
	\label{obj_dist}
\end{align}
where $\beta$ is a weight coefficient. Since the distance EVs can move during a given time interval is limited, we have the following constraints for variables $X^k$ and $Y^k$.
\begin{align}
	\begin{split}
		x_{ij}^k \geq 0 \text{ and } x_{ij}^k = 0 \text{ when } w_{ij} \geq m_1;\\
        y_{ij}^k \geq 0 \text{ and } y_{ij}^k = 0 \text{ when } w^*_{ij} \geq m_2,
%	    X^k \circ M^k =0,\quad X^k_{ij} \geqslant 0 %\quad i,  j \in \{1, 2, \dots, n^k\}  
	\end{split}
	\label{bound}
\end{align}
where $m_1 > 0, m_2>0$ is the upper bound balancing distance for an vacant and low-battery EV, respectively. These constraints consider real-world scenarios that we cannot dispatch EV to some far away regions when the moving distance exceeds the capability, either due to speed limit or insufficient battery.

\subsection{Utilization of Charging Stations $J_E$}
Given fixed number and locations of charging stations in the city, to avoid the long waiting time at some charging stations, one method is to balance the charging station utilization across the whole city. It will improve charging efficiency of EVs as well as decrease driver's potential cost due to EV's unique charging problems. 
The supply variable $c^k_i$ is the service rate or the average number of new available charging spots (due to EVs finished charging) in region $i$ during time $k$. The net total number of low-battery EVs $Y_i^k = \sum\limits_{j=1}^{N} y^k_{ji}-\sum\limits_{j=1}^{N} y^k_{ij}$ in region $i$ after balancing according to decision variable $Y^k$ is the average arrival rate or the average number of arriving EVs. The overall utilization of charging stations in region $i$ during time $k$ is approximated as $\frac{Y_i^k}{c_i^k}$. The total difference between the inverse local and inverse global utilization for $\tau$ time intervals is
\begin{align}
	\sum_{k=1}^{\tau}\sum_{i}^{N}\left| \frac{c^k_i}{Y^k_i} - \frac{\sum_{j=1}^{N}c^k_j}{\sum_{j=1}^{N} Y^k_i}\right|.
	\label{quality_charging_diff}
\end{align}
However, the charging station supply is a random variable that the prediction error can not be ignored~\cite{AMoD_queue, Tian2016Real, wang2019shared}, and we define the uncertainty set as $c\sim F^*_c, F^*_c\in \mathcal{F}_c$. Function~\eqref{quality_charging_diff} is not concave over uncertainty parameter $c^k$, for computationally tractability, we consider minimizing the following utilization quality function $J_E$
\begin{align}
	\begin{split}
		J_E(Y^{1:\tau}) = \sum\limits_{k=1}^{\tau} \sum\limits_{i=1}^{N} \frac{c_i^k}{(Y_i^k)^a}.
	\end{split}
	\label{quality_charging_one}
\end{align}
According to Lemma 1 in~\cite{ddrobust_Miao}, when the power parameter $a>0$ is designed to be small enough, objective function~\eqref{quality_charging_one} is linear in $c^k$ and convex in $Y^{1:\tau}$ and approximates the objective~\eqref{quality_charging_diff}.

\subsection{Constraints Definitions}%
%An EV is "available to be dispatched" means this EV is vacant and has enough power to begin a next trip.
\subsubsection{Service quality metrics}\label{sec:Service} Demand-supply ratio is one service quality metric for AMoD systems~\cite{mpcmod_icra16, ddrobust_Miao}. In this work we also minimize the total difference between local and global demand-supply ratio for $\tau$ time intervals. Let $V_i^k, O_i^k \in \mathbb{R}_+$ be the number of vacant and occupied EVs respectively at region $i$ at the beginning of time $k$ before balancing, and $V^k, O^k \in \mathbb{R}_+^N$. Define $S_i^k > 0$ as the total number of supply EVs available to be dispatched in region $j$ during time $k$, dispatch decisions as $X^{1:k}=\{X^1,X^2,...,X^k\}$, $Y^{1:k}=\{Y^1,Y^2,...,Y^k\}$. Then the following equations of $V_i^k, O_i^k, S_i^k$ describe dynamics for $\tau$ time steps:
\begin{align}
	\begin{split}
		 S^k_i &=\sum\limits_{j=1}^{N} x^k_{ji}-\sum\limits_{j=1}^{N} x^k_{ij} + V^k_i = X_i^k + V^k_i >0,\\ %k&=1,\dots,\tau, \\
		V^{k+1}_i&= \sum\limits_{j=1}^{N} P^k_{vji}S^k_j+ \sum\limits_{j=1}^{N} Q^k_{vji}O^k_j + c_i^k,\\ %k=1,\dots,\tau-1, \quad
		 O^{k+1}_i&=\sum\limits_{j=1}^{N}P^k_{oji}S^k_j+ \sum\limits_{j=1}^{N} Q^k_{oji}O^k_j,%\\ k&=1,\dots,\tau-1,  
	\end{split}
	\label{trans}
\end{align}
where $X_i^k = \sum\limits_{j=1}^{N} x^k_{ji}-\sum\limits_{j=1}^{N} x^k_{ij}$ is the net change of available EVs due to decision variable $X^k$ at region $i$. $P^k_v, P^k_o, Q^k_v, Q^k_o \in \mathbb{R}^{N \times N}$ are region transition matrices: $P^k_{vji}(P^k_{oji})$ is the probability that a vacant EV moves from region $j$ at the beginning of time $k$ will transverse to region $i$ and being vacant (occupied) at the beginning of time $k+1$, respectively. Similarly, $Q^k_{vji}$ and $Q^k_{oji}$ are the probability that an occupied EV moves from region $j$ at time $k$ will go to region $i$ and being vacant and occupied at the beginning of time $k+1$, respectively. When receding the time horizon, GPS locations (region information) and status of all EVs will always be updated by real-time sensing data and $V^1, O^1$ are provided by real-time data. We consider the following service quality constraints to make sure the demand-supply ratio of each region is within a similar range for service fairness:
\begin{align}
	\begin{split}
       l_i^k \leq \frac{r_i^k}{S_i^k} \leq h_i^k, \quad k=1,\dots,\tau,
	\end{split}
	\label{quality_region}
\end{align}
where $l_i^k(h_i^k)$ is the lower(upper) bound of the supply-demand ratio in region $i$ at time $k$. The value of $l_i^k$ and $h_i^k$ are decided by historical data. We transfer the inequalities~\eqref{quality_region} to the following equations form with slack variables $D_i^k, U_i^k$:  
\begin{align}
	\begin{split}
		r_i^k - l_i^k S_i^k - (D_i^k)^2 &= 0,\\
        r_i^k - h_i^k S_i^k + (U_i^k)^2 &= 0, \quad k=1,\dots,\tau.
	\end{split}
	\label{quality_equ}
\end{align}

\iffalse
So the objective would be minimizing the total difference between the local region demand-supply ratio to the system-level demand-supply ratio:
\begin{align}
	\sum_{k=1}^{\tau}\sum_{i}^{N}\left| \frac{r_i^k}{S^k_i} - \frac{\sum_{j=1}^{N}r^k_j}{\sum_{j=1}^{N} S^k_j}\right|.
	\label{quality}
\end{align}
However, function~\eqref{quality} is not concave over the uncertainty parameter $r^k$, and one minmax problem with~\eqref{quality} as part of the objective function is computationally intractable when the max is over $r^k$~\cite{Ye_dro, ddrobust_Miao}.
\fi

\subsubsection {Constraints on decision variables $Y^{1:\tau}$} We define $L_i^k \in \mathbb{R}_+$ be the total amount of low-battery EVs in region $i$ before balancing at the beginning of time $k$. Then $L_i^k$ should have the following relationship with $S_i^k$ and 
$Y^k$ %for $k=1,\dots,\tau-1$
\begin{align}
	\begin{split}
		L^{k+1}_i &=\sum\limits_{j=1}^{N} y^k_{ji}-\sum\limits_{j=1}^{N} y^k_{ij} + \sum\limits_{j=1}^{N} P_{lji}^k S^k_{j} >0,
	\end{split}
	\label{con_charing}
\end{align}
where $L^k \in \mathbb{R}_+^N$ and $L^1$ is given by real-time data, $P^k_l \in \mathbb{R}^{N \times N}$ is the region transition matrix: $P^k_{lji}$ is the probability that a vacant EV moves from region $j$ at the beginning of time $k$ will go to region $i$ and being low-battery status at the beginning of time $k+1$. Thus the region transition matrices estimated from data satisfy that $\sum\limits_{j=1}^{N} P^k_{lij} + P^k_{vij} + P^k_{oij} = 1$ and $\sum\limits_{j=1}^{N} Q^k_{oij} + Q^k_{vij} = 1$. %The method of calculating these matrices are introduced in~\cite{taxi_tase16}. 

\subsection{Predicted Model and Uncertainty Set Construction}
\label{pmus}
Instead of assuming we know the true probability distributions of $r$ and $c$ from data, in this work we construct uncertainty sets $\mathcal{F}_r$ and $\mathcal{F}_c$ that describe possible probability distributions of $r$ and $c$ by applying the Algorithm 2 proposed in~\cite{dddro_tcps20} which constructs distributional sets with a general prediction model. Here we use the autoregressive integrated moving average model (ARIMA model)~\cite{ARIMA} as the prediction model to capture the spatial and temporal correlations for the predictions. The coefficients of ARIMA model can be estimated by maximum likelihood estimation. We use $\hat{r}$ and $\hat{c}$ to denote the the predicted value, $\tilde{r}$ and $\tilde{c}$ to denote the sample value or $r$ and $c$, respectively. And then we get a corresponding estimation residuals as
\begin{align}
	\begin{split}
	    \tilde{\delta}_r = \tilde{r} - \hat{r}, \quad \tilde{\delta}_c = \tilde{c} - \hat{c}.
	\end{split}
	\label{residual}
\end{align}

Then we use these estimation residuals to construct distribution sets $\mathcal{F}_r$ and $\mathcal{F}_c$ that contain the true distribution of $r$ and $c$ with probability at least $1-\alpha_h$ where $\alpha_h$ is a significant values can be freely chosen from 0 to 1. $\mathcal{F}_r$ and $\mathcal{F}_c$ have the following format:

\footnotesize{
\begin{align}
\begin{split}
\mathcal{F}_r(\hat{r}, \hat{\Sigma}_r, \hat\gamma_{1r}, \hat\gamma_{2r}) = \{ &r=\hat{r} + {\delta}_r: (\mathbb{E}[\delta_r])^T \hat{\Sigma}_{r}^{-1}\mathbb{E}[\delta_r] \leqslant \hat\gamma_{1r},\\ 
&\mathbb{E}(\delta_r \delta_r^T) \leq \hat\gamma_{2r} \hat{\Sigma}_{r}\};
\\
\mathcal{F}_c(\hat{c}, \hat{\Sigma}_c, \hat\gamma_{1c}, \hat\gamma_{2c}) = \{ &c=\hat{c} + {\delta}_c: (\mathbb{E}[\delta_c])^T \hat{\Sigma}_{c}^{-1}\mathbb{E}[\delta_c] \leqslant \hat\gamma_{1c}, \\
&\mathbb{E}(\delta_c \delta_c^T) \leq \hat\gamma_{2c} \hat{\Sigma}_{c}\},
\end{split}
\label{uncertain_delta}
\end{align}}
\normalsize
where $\delta_r (\delta_c)$ is the difference between true value $r(c)$ and predicted value $\hat{r} (\hat{c})$. Since we don't have true value of the difference $\delta_r = {r} - \hat{r} \ (\delta_c = {c} - \hat{c})$, we use estimation residuals to capture the information from these difference. $\hat{\Sigma}_r(\hat{\Sigma}_c)$ is the estimated covariance, $\hat\gamma_{1r}, \hat\gamma_{2r} \ (\hat\gamma_{1c}, \hat\gamma_{2c})$ are two estimated threshold values of $\delta_r(\delta_c)$ by Algorithm 2 in~\cite{dddro_tcps20} based on the concept of bootstrapping. More discussions about uncertainty set construction refer to~\cite{Ye_dro, dddro_tcps20}.

\subsection{Distributionally Robust EV Balancing Problem}
The final goal is to dispatch EVs under minimal cost. We define a weight parameter $\theta$ of the two objectives $J_D$ defined in \eqref{obj_dist} and $J_E$ defined in~\eqref{quality_charging_one}. With constraints~\eqref{bound}, \eqref{trans}, \eqref{quality_equ}, \eqref{con_charing}, we define the following distributionally robust EVs balancing problem under uncertain probability distributions of random demand and supply:
\begin{align}
	\begin{split} 
		\min_{\substack{X^{1:\tau}, Y^{1\tau}, S^{1:\tau}, D^{1:\tau};\\ U^{1:\tau},  V^{2:\tau}, O^{2:\tau}, L^{2:\tau}}}
		\underset{\{F_r\in \mathcal{F}_r,F_c\in \mathcal{F}_c\} }{\text{max}}\ &\mathbb{E}\left[J_D+\theta J_E \right]\\
		\text{s.t.}\quad
		%& X^k_{ij} W^k_{ij} \leq m^k X^k_{ij}, \\
		\text{\eqref{bound}, \eqref{trans}, \eqref{quality_equ}, \eqref{con_charing}}.
	\end{split}
	\label{obj_final}
\end{align}
Since above problem~\eqref{obj_final} can not be calculated in polynomial time directly, we derive a computationally tractable form of this problem in the following section.

\section{Computationally Tractable Form}
\label{sec:algorithm}
In this section, we derive the theoretical result, Theorem 1 of this work, i.e., a computationally tractable and equivalent convex optimization form for problem~\eqref{obj_final} via strong duality. 
Hence, the optimal solution of~\eqref{obj_final} can be calculated in real time considering both passenger demand and EV supply uncertainties.

\begin{theorem}
	The distributionally robust optimization problem~\eqref{obj_final} with two distributional sets~\eqref{uncertain_delta} is equivalent to the following convex optimization problem
	\begin{align}
		\begin{split}
			&\min_{\substack{X^{1:\tau}, Y^{1:\tau}, D^{1:\tau}, U^{1:\tau};\\ S^{1:\tau}, V^{2:\tau}, O^{2:\tau},L^{2:\tau}; \\ {Q_r,q_r,v_r,t_r,Q_c,q_c,v_c,t_c}}}
			\quad  H_o + v_r + t_r + v_c + t_c\\
			&\text{s.t.}
			\quad  \begin{bmatrix}v_r & \frac{1}{2}(q_r+\lambda_U + \lambda_D)^T\\ \frac{1}{2}(q_r+\lambda_U + \lambda_D) & Q_r 
			\end{bmatrix} \succeq 0,\\
			&\quad \quad \begin{bmatrix}v_c & \frac{1}{2}(q_c+\lambda_V - Z)^T\\ \frac{1}{2}(q_c+\lambda_V - Z) & Q_c 
			\end{bmatrix} \succeq 0,\\
			&\quad \quad t_r \geqslant (\hat\gamma_{2r} \hat{\Sigma}_r+\hat{r}\hat{r}^T)\cdot Q_r + \hat{r}^T q_r\\ &\quad \quad \quad +\sqrt{\hat\gamma_{1r}} \|\hat{\Sigma}_r^{1/2} (q_r+2Q_r\hat{r})\|_2 ,\\
			&\quad \quad t_c \geqslant (\hat\gamma_{2c} \hat{\Sigma}_c+\hat{c}\hat{c}^T)\cdot Q_c + \hat{c}^T q_c\\ &\quad \quad \quad +\sqrt{\hat\gamma_{1c}} \|\hat{\Sigma}_c^{1/2} (q_c+2Q_c\hat{c})\|_2 ,\\
			&\quad \quad Q_r, Q_c, \lambda, v_r, v_c \succeq 0 ,\\
			&\quad \quad Z_i^k \geqslant \frac{1}{(Y_i^k)^a}\\
			&\quad \quad x_{ij}^k \geqslant 0 \text{ and } x_{ij}^k = 0 \text{ when } w_{ij} \geq m;\\
            &\quad \quad y_{ij}^k \geqslant 0 \text{ and } y_{ij}^k = 0 \text{ when } w^*_{ij} \geq m,
		\end{split}
		\label{thm1}
	\end{align}
	where $H_o=J_D -  ( \lambda_S^T f_S + \lambda_O^T f_O + \lambda_L^T f_L + \lambda_s^T S^{1:\tau} + \lambda_l^T L^{2:\tau}) - \lambda_D^T diag(lS^T - DD^T) - \lambda_U^T diag(hS^T + UU^T) - \sum\limits_{k = 1}^{\tau -1}(-V^{k+1}_i + \sum\limits_{j=1}^{N} P^k_{vji}S^k_j+ \sum\limits_{j=1}^{N} Q^k_{vji}O^k_j)\lambda_{V_i^k} $,
	$J_D$ is defined as~\eqref{obj_dist}.
	\label{theorem1}
\end{theorem}
\textbf{Proof}. See Appendix~\ref{proof_theorem1}.
%To prove Theorem~\ref{theorem1}

\section{Evaluation with Real-World E-taxi Data}
\label{sec:simulation}
In this section, we evaluate the performance of the proposed distributionally robust optimization-based EV balancing algorithm~\eqref{theorem1} with three-week E-taxi data from the Chinese city Shenzhen (one of the largest cities in China that operates over 10,000 E-taxis). 

In total, there are four different datasets used in this paper, including E-taxi GPS data (vehicle ID, locations, time and speed, etc), transaction data (vehicle ID, pick-up and drop-off time, pick-up and drop-off location, travel distance, etc), charging station data (locations, name, the number of charging points, etc), and the urban partition data (geographic boundaries of 491 small separate regions composing Shenzhen).

We first utilize a widely adopted spatiotemporal constraint-based method \cite{Li2015Growing, wang2019shared, Tian2016Real} to extract charging events of E-taxis and re-split Shenzhen into 54 large areas based on these data. After obtaining the charging events, we utilize the first two weeks as the training data to determine the uncertainty sets and parameters. The remaining one week is used as the testing data to compare the vehicle balancing cost. 

\subsection{Uncertainty Set}
Here we select 5 busy regions' data as the input for Algorithm 2 in~\cite{dddro_tcps20} to construct uncertainty sets. 
We set time horizon $\tau$ as $2$ and a significant $\alpha_h = 0.25$ to make sure the true probability distributions of demand $r$ and supply $c$ are separately contained in the constructed uncertainty sets with probability at lease $75\%$. We show how thresholds $\hat\gamma_{1c}$, $\hat\gamma_{2c}$ change for different sample in table~\ref{tab:gamma}. The value of $\hat\gamma_{1c}$ and $\hat\gamma_{2c}$ decrease when $N_B$ increase and the speed of decreasing becomes slow as $N_B$ becoming larger and larger. We also notice that as $N_B$ increase, the corresponding estimated moments have a trend to converge to a certain constant which meets bootstrapping algorithm's intuition. 
% Fig.~\ref{fig:norm_dis} shows the centralized supply uncertainty set is convex and closed and the board of this set follows a multivariate Gaussian distribution. The contour lines in Fig.~\ref{fig:norm_dis} shows the board of constructed uncertainty set while the color indicates the numerical value of the board: darker the color, higher the value. 

\begin{table}[]\centering
\vspace*{0pt}
\caption{Thresholds $\hat\gamma_{1c}$ and $\hat\gamma_{2c}$ for different sample number $N_B$}
\vspace*{-6pt}
\begin{tabular}{|c|c|c|c|c|c|c|}
\hline
$N_B$     & 10    & 20    & 50    & 100   & 500   & 1000  \\ \hline
$\hat\gamma_{1c}$ & 1.504 & 0.964 & 0.576 & 0.399 & 0.296 & 0.176 \\ \hline
$\hat\gamma_{2c}$ & 3.715 & 2.832 & 2.006 & 1.768 & 1.374 & 1.317 \\ \hline
\end{tabular}
\label{tab:gamma}
\end{table}

\subsection{Predicted Model}

\begin{figure}[!t]
	\centering
	\includegraphics [width=0.42\textwidth]{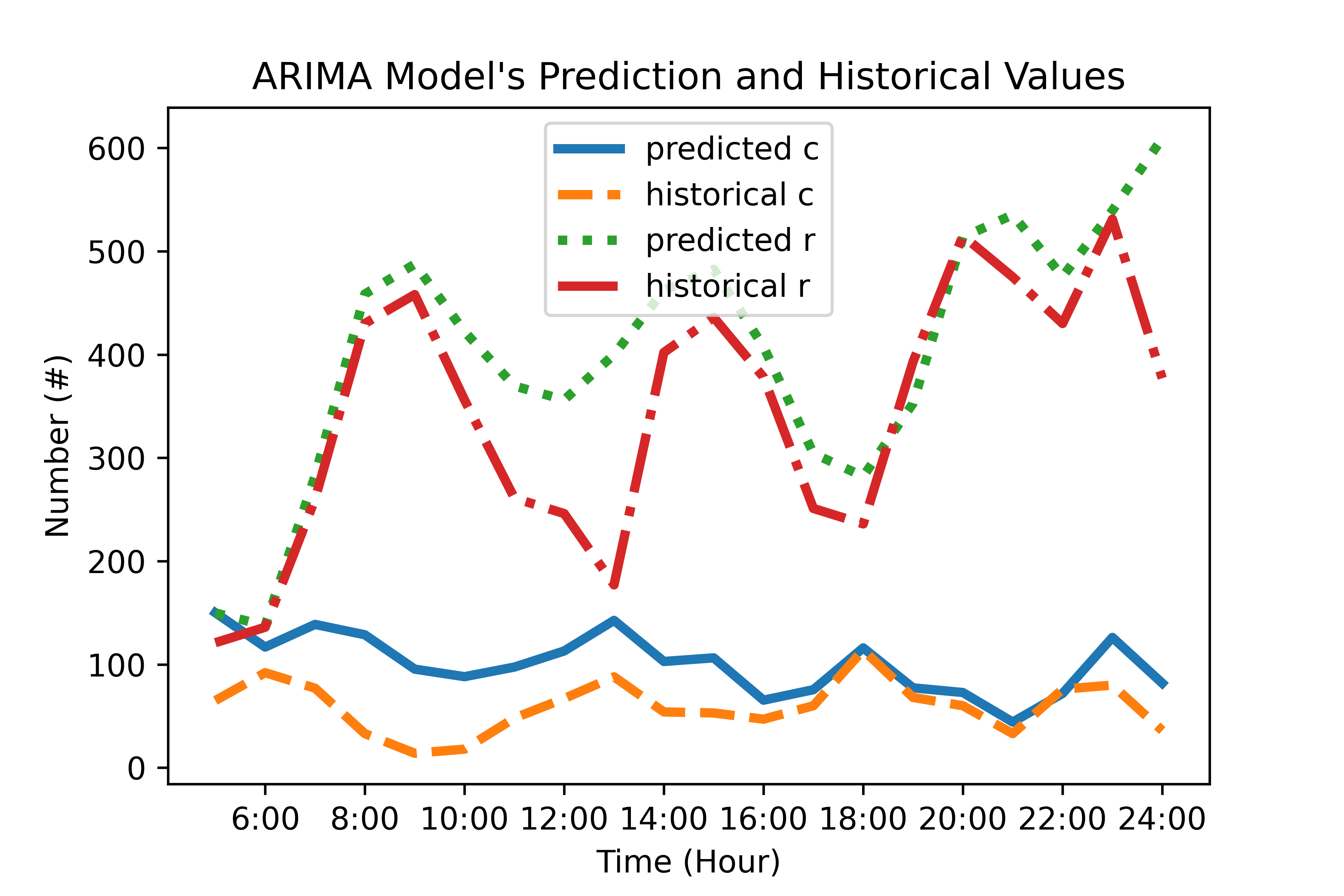}
	\vspace{-10pt}
	\caption{ARIMA model demonstrates time trends very well.}
	\label{fig:arima}
	\vspace{-8pt}
\end{figure}

As mentioned in section~\ref{pmus}, we use ARIMA model~\cite{ARIMA} to predict $r$ and $c$ where current values are predicted by former concatenation of values. Fig.~\ref{fig:arima} compares ARIMA model's predicted values and historical values of $c$ and $r$ on one day. We can see that ARIMA model demonstrates data's time trends very well and the predicted values of $r$ in peak hours: 8am-10am, 2pm-4pm, 8pm-10pm, are very close to historical values.

\subsection{The Performance of the Proposed Method}
\begin{figure}[!t]
	\centering
	\includegraphics [width=0.42\textwidth]{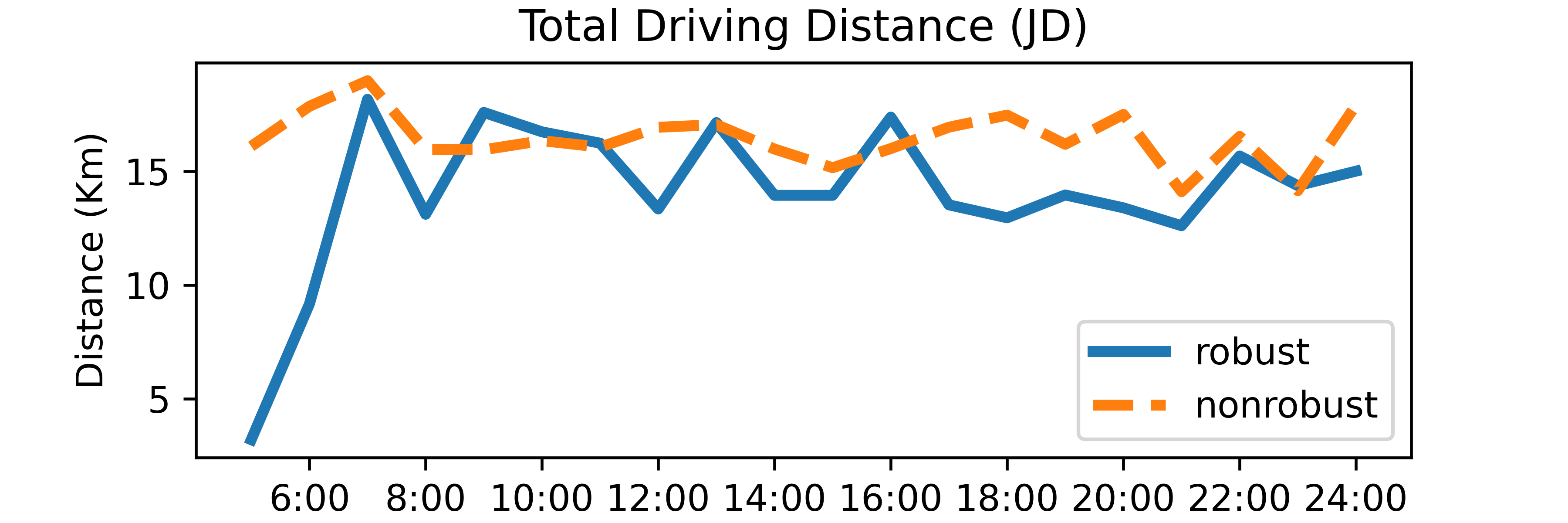}
	\vspace{-10pt}
	\caption{By using robust optimization, the average total driving distance is reduced by 14.49\% compared to NonRobust method.}
	\label{fig:total_cost}
	\vspace{-8pt}
\end{figure}

We compare our distributionally robust method with the non-robust method~\cite{p2charge} by using the same real-time sensing data to evaluate the performance. In the robust method we use the demand and supply uncertainty sets defined as~\eqref{uncertain_delta}, while in the non-robust method, they are deterministic vectors. In Fig~\ref{fig:total_cost}, we compare the total driving distance of applying the optimal decision of each method from 5am to next day's 12am, 19 hours in total. Here, the total driving distance $J_D$ in~\eqref{obj_dist} is defined as a weighted sum of the charging idle distance and the service idle distance for one EV. Total driving distance is supposed to be small since it implies low total cost. We can see that most of time, the total cost of the robust model is lower than that of the non-robust model. In particular, the average total driving distance is reduced by 14.49\% compared with the non-robust method. 

\begin{figure}[!t]
	\centering
	\includegraphics [width=0.42\textwidth]{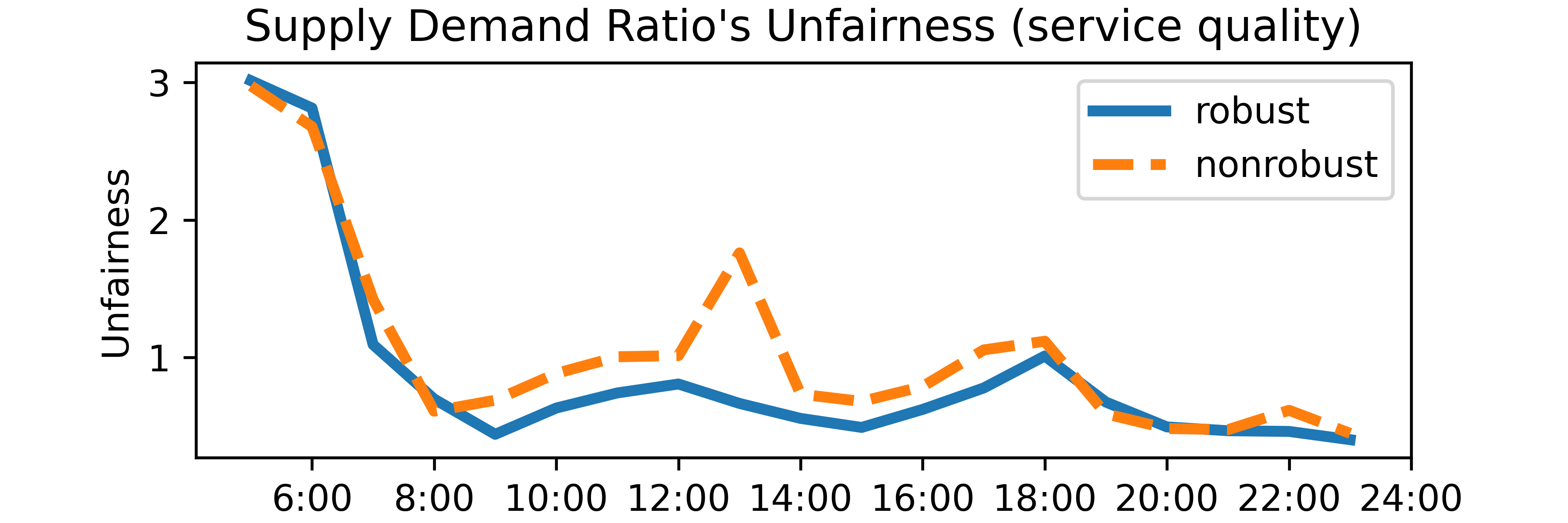}
	\vspace{-10pt}
	\caption{By using robust optimization, the average unfairness of supply-demand ratio is reduced by 15.78\% compared to NonRobust method.}
	\label{fig:fairm}
	\vspace{-8pt}
\end{figure}

\begin{figure}[!t]
	\centering
	\includegraphics [width=0.42\textwidth]{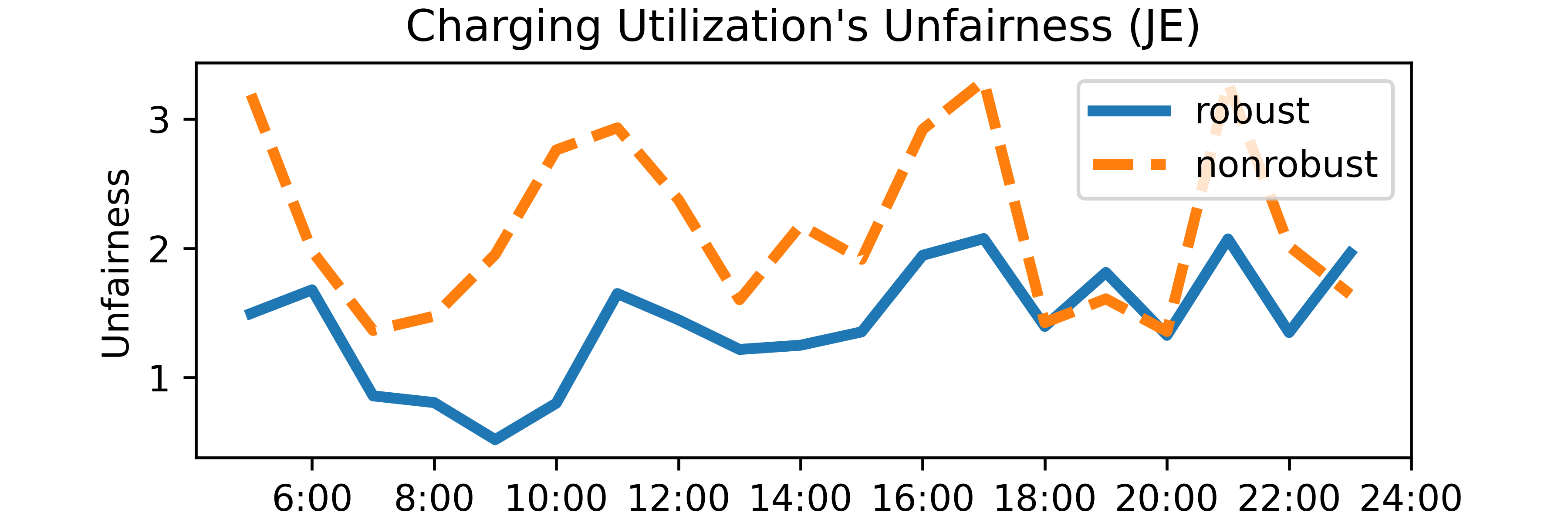}
	\vspace{-10pt}
	\caption{By using robust optimization, the average unfairness of utilization is reduced by 34.51\% compared to NonRobust method.}
	\label{fig:fairc}
	\vspace{-8pt}
\end{figure}

In Fig~\ref{fig:fairm} and~\ref{fig:fairc}, we compare the unfairness of supply-demand ratio and utilization of the whole city taking decisions of each method. The unfairness metric of utilization (supply-demand ratio) is designed as the total sum of absolute difference between the inverse-local utilization (local supply-demand ratio) and inverse-global utilization (global supply-demand ratio). Lower sum of absolute difference means higher fairness but lower unfairness. The unfairness of utilization in~\eqref{quality_charging_diff} is actually accurate version of $J_E$ in~\eqref{quality_charging_one} before we approximate $J_E$. The unfairness of supply-demand ratio in section~\ref{sec:Service} is the service quality metric for taxi dispatch in~\cite{Miao_iccps17}: the lower the service quality metric, the better customers' accumulative satisfaction is since they can enjoy uniform quality service. By using robust optimization, the average unfairness of supply-demand ratio and utilization is reduced by 15.78\% and 34.51\% respectively compared to NonRobust method.

%\subsection{}

\section{Conclusion}
\label{sec:conclusion}
Autonomous mobility-on-demand systems can provide more efficient services, and the total idle distance can be reduced with vehicle balancing algorithms in general. However, with the increasing amount of EVs and the limited charging facilities in the city, the uncertainty of charging time (including waiting time) at a charging station affects the EVs supply to provide efficient service for AMoD systems. In this paper, we design a data-driven distributionally robust EV balancing method to minimize the worst-case expected cost under uncertainties about the probability distributions of both demand and supply. Besides reducing EVs' total idle driving distance, we also balance the demand-supply ratios and the charging station utilization of different regions among the city. Then we prove an equivalent computationally tractable form of the distributionally robust problem under the ellipsoid uncertainty sets constructed from data. Evaluations based on three-week real-world E-taxi data from the Chinese city Shenzhen show that the average total balancing cost is reduced by 14.49\%, the  average  unfairness  of  supply-demand ratio and utilization is reduced by 15.78\% and 34.51\%, respectively. In the future, we will further evaluate our algorithm based on large-scale data of several years from multiple cities.

\bibliographystyle{abbrv}
{ \small 
\bibliography{ref}
}

\section{Appendix}
\subsection{Proof of Theorem~\ref{theorem1}}
\label{proof_theorem1}
\begin{proof}
We notice that the uncertainty parameters are involved in both objective functions and constraints. We first use $\underset{F_r\in \mathcal{F}_r}{\text{max.}} \mathbb{E}(r_i^k)$ substitute $r_i^k$, $\underset{F_c\in \mathcal{F}_c}{\text{min.}} \mathbb{E}(c_i^k)$ substitute $c_i^k$ in constraints \eqref{trans} and \eqref{quality_equ}. It's reasonable due to the idea of minimizing the worst case: any uncertain values of $r_i^k$ and $c_i^k$ should meet the relationship with other decision variables shown in constraints \eqref{trans} and \eqref{quality_equ}. No matter which probability distribution is selected as the specific distribution to attain the worst case, the simplest worst case for single value of $r_i^k$ and $c_i^k$ is the case that the demand is really large that attains the maximal possible demand value while the supply is very small that attains the minimal possible supply value. And $\mathbb{E}(r_i^k) (\mathbb{E}(c_i^k))$ is the probability-weighted average of all its possible values. Then we transfer all constraints into functional formats as below:
    \footnotesize{
    \begin{align}
	\begin{split}
	    f_{D_i^k} &= \underset{F_r\in \mathcal{F}_r}{\text{max.}} \mathbb{E}(r_i^k) - l_i^k S_i^k - (D_i^k)^2 = 0,\\
        f_{U_i^k} &= \underset{F_r\in \mathcal{F}_r}{\text{max.}} \mathbb{E}(r_i^k) - h_i^k S_i^k + (U_i^k)^2 = 0,\\ 
	    f_{S_i^k} &= -S^k_i + X_i^k + V^k_i = 0,\quad k=1,\dots,\tau, \\
		f_{V_i^{k+1}} &= -V^{k+1}_i + \sum\limits_{j=1}^{N} P^k_{vji}S^k_j+ \sum\limits_{j=1}^{N} Q^k_{vji}O^k_j \\
		& \quad - \underset{F_c\in \mathcal{F}_c}{\text{max.}} \mathbb{E}(-c_i^k) = 0,\\ 
		f_{O_i^{k+1}} &= -O^{k+1}_i + \sum\limits_{j=1}^{N}P^k_{oji}S^k_j+ \sum\limits_{j=1}^{N} Q^k_{oji}O^k_j = 0,,\\
	    f_{L_i^{k+1}} &= - L^{k+1}_i + Y_i^k + \sum\limits_{j=1}^{N} P_{lji}^k S^k_{j} = 0,\\ %\quad k=1,\dots,\tau-1,\\
	     L_i^k &> 0, \quad k=1,\dots,\tau-1,\quad
	     S_i^k > 0, \quad k=1,\dots,\tau.
	\end{split}
	\label{con_functional}
\end{align}}
\normalsize
	Let $f_D = [f_{D_1^1},f_{D_1^2},...,f_{D_1^\tau},...,f_{D_N^\tau}]^T \in \mathbb{R}^{N\tau}$ be a constraint function vector, for $i = 1, \dots,N, k = 1,\dots,\tau-1$, and $f_U,f_S,f_V,f_O,f_L$ have the same definition but for computational convenient, if one's dimension is less than $N\tau$, we add 0 in corresponding missing positions to complete its dimension. For the primal maximization problem $\underset{F_r\in \mathcal{F}_r,F_c\in \mathcal{F}_c }{\text{max}}\ \mathbb{E}\left[J_D+\theta J_E \right],\quad
		\text{s.t } \eqref{con_functional},$\\
%	\begin{align*}
%	\end{align*}
	its associated Lagrange dual problem \eqref{dual} can obtain its best upper bound since the strong duality holds. The primal objective function only contains uncertainty parameter $c$ and is concave over $c$ because it's a linear function of $c$ when given other decision variables. The constraints are also all linear in $r$ and $c$. When the primal problem is in this case, we usually have strong duality, and Slater’s theorem also states that strong duality holds according to~\cite{Bookcvx_Boyd}. More context about duality are discussed in~\cite{Bookcvx_Boyd}.
	\begin{align}
	\begin{split}
	   & \underset{\lambda \succeq 0}{\text{min}}\
		\underset{F_r\in \mathcal{F}_r,F_c\in \mathcal{F}_c }{\text{max}}\ J_{dual},\\
%		\end{split}
%	\end{align}
%	where
%	\begin{align}
%	\begin{split}
	    J_{dual} &= \mathbb{E}\left[J_D+\theta J_E \right] - (\lambda_U^T f_U + \lambda_S^T f_S + \lambda_V^T f_V \\
	    &+ \lambda_O^T f_O + \lambda_L^T f_L + \lambda_s^T S^{1:\tau} + \lambda_l^T L^{2:\tau})
	\end{split} 
	\label{dual}
	\end{align}
	$\lambda_U^T, \lambda_S^T, \lambda_V^T, \lambda_O^T, \lambda_L^T, \lambda_s^T, \lambda_l^T$ are corresponding Lagrange multipliers and $\lambda$ is defined as a vector combined by all these Lagrange multipliers. We have $\frac{c_i^k}{(Y^k_i)^a} \geqslant 0$ and $c_i^k \geqslant 0$ by the definitions of $J_E$ in~\eqref{quality_charging_one} and the Queuing model, then for any vector $Z\in \mathbb{R}^{N\tau}$, $Z=[z^1_1, z^1_2,\dots, z^{\tau}_1, z^{\tau}_2, \dots,z^{\tau}_{N\tau}]^T$ that satisfies $0<\frac{1}{(Y^k_i)^{a}} \leqslant z_{i}^k$, we also have\\
	\centerline{$
		%	\begin{align*}
		0 \leqslant \sum_{k=1}^{\tau} \sum\limits_{i=1}^{N} \frac{c^k_i}{(Y^k_i)^a} \leqslant Z^T c,
		%	\end{align*}
		$}
	and the second inequality strictly holds when all $\frac{r^k_i}{(Y^k_i)^{a}}=z_i^k$, for $i=1,\dots,N$, $k=1,\dots, \tau$. The constraints of problem~\eqref{dual} are independent of $c$, hence, for any $c$, the minmax problem \eqref{dual} is equivalent to 
	\begin{align}
	    \begin{split}
	        \underset{\lambda \succeq 0}{\text{min}}\
		    \underset{F_r\in \mathcal{F}_r,F_c\in \mathcal{F}_c }{\text{max}}\ &J_{dual}^\prime\\
		    &\text{s.t.} \frac{1}{(Y_i^k)^a} \leqslant z_i^k, \quad Z \in \mathbb{R}^{N\tau}
		\end{split}
	    \label{dual_2}
	\end{align}{}
	where
	\begin{align}
	\begin{split}
	    J_{dual}^\prime &= \mathbb{E}\left[J_D+\theta Z^T c \right] - (\lambda_U^T f_U + \lambda_S^T f_S + \lambda_V^T f_V \\
	    &+ \lambda_O^T f_O + \lambda_L^T f_L + \lambda_s^T S^{1:\tau} + \lambda_l^T L^{2:\tau})
	\end{split}
	\end{align}
	In the dual problem's objective function $J_{dual}^\prime$, not all parts contain uncertainty parameters. We can separate $J_{dual}^\prime$ into three parts
	\begin{align}
	\begin{split}
	  	H_r &=  -(\lambda_U^T + \lambda_D^T) r, \quad
	    H_c =  \theta J_E - \lambda_V^T c, \\
	    H_o &= J_{dual}^\prime - \mathbb{E}[H_c + H_r]  ,
	\end{split}
	\end{align}
	Where only $H_r$ contains all $r$, $H_c$ contains all $c$. $H_o$ can be put as a deterministic value when given other decision variables. So we can turn to consider the following maximization problem
	\begin{align}
		\underset{r \sim F_r, c \sim F_c, F_r\in \mathcal{F}_r,F_c\in \mathcal{F}_c}{\text{max}} \mathbb{E} [H_r + H_c].
		\label{obj_max_all}
	\end{align}
	Since $r$ and $c$ are independent, problem \eqref{obj_max_all} equals the separated maximization problem
	\begin{align}
		\underset{r \sim F_r, F_r\in \mathcal{F}_r}{\text{max}} \mathbb{E} [H_r] + \underset{c \sim F_c, F_c\in \mathcal{F}_c}{\text{max}} \mathbb{E}[H_c].
		\label{obj_max_sep}
	\end{align}
	Problem \eqref{obj_max_sep} satisfies the conditions of Lemma 1 in~\cite{Ye_dro}, and the maximum expectation value of $H_r + H_c$ for any possible $r \sim F_r ,c \sim F_c$ where $F_r \in \mathcal{F}_r ,F_c \in \mathcal{F}_c$ equals the optimal value of the problem
		\begin{align}
		\begin{split}
			\min_{\substack{Q_r,q_r,v_r,t_r;\\Q_c,q_c,v_c,t_c}}\quad &v_r + t_r + v_c + t_c\\
			\text{s.t.}
			\quad & v_r\geqslant H_r -r^T Q_r r -r^T q_r,\\
			&t_r \geqslant (\hat\gamma_{2r} \hat{\Sigma}_r+\hat{r}\hat{r}^T)\cdot Q_r + \hat{r}^T q_r \\&\quad\quad+\sqrt{\hat\gamma_{1r}} \|\hat{\Sigma}_r^{1/2} (q_r+2Q_r\hat{r})\|_2 ,\\
			\quad & v_c\geqslant H_c -c^T Q_c r -c^T q_c,\\
			&t_c \geqslant (\hat\gamma_{2c} \hat{\Sigma}_c+\hat{c}\hat{c}^T)\cdot Q_c + \hat{c}^T q_c \\&\quad\quad+\sqrt{\hat\gamma_{1c}} \|\hat{\Sigma}_c^{1/2} (q_c+2Q_c\hat{c})\|_2 ,\\
			&Q_r, Q_c\succeq 0 .\\
		\end{split}
		\label{step1}
	\end{align}
	Note that the first and third constraints about $v_r$ and $v_c$ is equivalent to $v_r \geqslant f_r(r^*)$ and $v_c \geqslant f_c(c^*)$ where $f_r(r^*) ( f_c(c^*) )$ is the optimal value of the following problem 
% 	(same form for $c$)
	    \begin{align}
	        \begin{split}
	            \max_{r} \quad & H_r - r^T Q_r r -r^T q_r\quad
	                 \text{s.t.} \quad r \geqslant 0,\\
	            \max_{c} \quad & H_c - c^T Q_c c -c^T q_c\quad
	                 \text{s.t.} \quad c \geqslant 0.
	        \end{split}{}
	        \label{obj_cons}
	    \end{align}{}
	Since $Q_r$ and $Q_c$ are positive semi-defined, $H_r (H_c)$ is a linear function over $r (c)$, problem \eqref{obj_cons} is convex. Solving this problem by taking partial derivative over $r (c)$ without constraints, we have:
	    \begin{align}
	        \begin{split}
	            v_r & \geqslant\frac{1}{4}(q_r + \lambda_U + \lambda_D)^T Q_r^{-1} (q_r + \lambda_U + \lambda_D)\\
	            v_c & \geqslant\frac{1}{4}(q_c + \lambda_V - Z)^T Q_c^{-1} (q_r + \lambda_V - Z)
	        \end{split}{}
	    \end{align}{}
	By Schur complement, the above constraints are
	    \begin{align}
	        \begin{split}
	            &\quad  \begin{bmatrix}v_r & \frac{1}{2}(q_r+\lambda_U + \lambda_D)^T\\ \frac{1}{2}(q_r+\lambda_U + \lambda_D) & Q_r 
			\end{bmatrix} \succeq 0,\\
			&\quad \quad \begin{bmatrix}v_c & \frac{1}{2}(q_c+\lambda_V - Z)^T\\ \frac{1}{2}(q_c+\lambda_V - Z) & Q_c 
			\end{bmatrix} \succeq 0,\\
	        \end{split}{}
	    \end{align}{}
	But this equivalent constraints are under the conditions of no constraints in problem \eqref{obj_cons}. We still have another two constraints for $v_r (v_c)$ that $v_r \geqslant f_r(r = 0) = 0$ and $v_c \geqslant 0$. Then use the fact that min-min operations can be performed jointly and combine all constraints we can reformulate problem \eqref{dual} as \eqref{thm1}.
\end{proof}

\end{document}

%% file: intro.tex
\section{Introduction}
There are over 5 million EVs by December 2018 globally, and this figure is predicted to increase to 100-125 million by 2030 \cite{EB2018}. Compared to conventional gas vehicles, EV fleets have prolonged charging time and concentrated mobility patterns due to current charging technologies and limited charging infrastructures, especially for commercial EV fleets, e.g., e-taxi, future autonomous mobility-on-demand (AMoD) systems, given their long daily travel distances \cite{epat}.

Researchers have been focusing on models and algorithms to study EVs~\cite{bCharge, epat}. There have also been focusing on how to choose the optimal locations for charging stations and how to assign charging points to EVs in each station to minimize the charging time of EVs considering various constraints, e.g., demand, costs, and charging compatibility~\cite{Li2015Growing}, \cite{wang2019shared}. However, high costs of charging infrastructures and land resources make it impractical to deploy abundant charging stations and points at the early promotion stage~\cite{wang2019shared}. Even when there is enough charging infrastructure for all EVs in theory, the uncontrolled and decentralized charging and mobility behaviors of some EV fleets, e.g., e-taxi, cause long waiting times when the demand for charging points greatly exceeds the availability~\cite{Tian2016Real}.

The above mentioned EV management issues have posed key optimization and scheduling algorithm challenges for world-wide EV adoption of AMoD. The interaction between AMoD systems and power networks through EVs based on the vehicles' charging requirements, battery depreciation, and power transmission constraints have been investigated, and the economic and societal value of EV AMoD has been analyzed~\cite{EVAMoD_tcns19}. To improve the performance of general AMoD systems, mobility demand based vehicle balancing methods have been proposed with various system design objectives, such as reducing the number of vehicles needed to serve all passengers~\cite{mpcmod_icra16, DDmpcmod_icra18}, reducing customers' waiting time~\cite{mod_iros18}, or taxis' total idle distance~\cite{taxi_tase16}. However, the limited knowledge we have about charging patterns~\cite{epat} affect the performance of vehicle balancing strategies, and make real-time decisions under demand model uncertainties still a challenging and unsolved task.

The contributions of this work are as follows:

\begin{itemize}
 	\item We are the first to consider both future demand uncertainties and EV supply uncertainties predicted based on charging activity data in designing a system-level vehicle balancing algorithm. While model predictive control algorithms~\cite{mpcmod_icra16, DDmpcmod_icra18} have been designed considering AMoD system demand uncertainties in the literature, the supply side uncertainties for EV AMoD is not well studied yet.
 	\item We design a distributionally robust optimization approach to balance EVs across a city for minimum total idle distance and balanced charging station utilization with respect to the worst-case expected cost. The approach considers probability distribution uncertainties of the passenger mobility demand and the EV supply caused by the challenge of charging process prediction~\cite{wang2019shared, Tian2016Real}. 
 	\item We derive an equivalent form of convex optimization problem for the proposed distributionally robust optimization problem to provide system-level performance guarantee in a computationally tractable way under model uncertainties. Based  on  real data of  Shenzhen  city,  we show  that the  average  total  balancing  cost is  reduced  by  14.49\%,  the  average  unfairness  of  supply-demand  ratio  and  utilization  is  reduced  by  15.78\%  and 34.51\%, respectively, with the proposed method, compared  with  solutions  which  do  not consider  model  uncertainties.
 \end{itemize}

The rest of the paper is organized as follows. The distributionally robust EV balancing problem is presented in Section~\ref{sec:prob_form}. An equivalent computationally tractable form is derived in Section~\ref{sec:algorithm}. We show performance improvement in experiments based on real data in Section~\ref{sec:simulation}. Concluding remarks are provided in Section~\ref{sec:conclusion}.